\documentclass[12pt,draftcls,onecolumn]{IEEEtran}
\usepackage{color}
\usepackage{ifpdf}
\ifpdf
  \usepackage[pdftex]{graphics}
\else
  \usepackage[dvips]{graphicx}
\fi
\IEEEoverridecommandlockouts

\usepackage{stmaryrd}
\usepackage{amsfonts}
\usepackage{amsmath}
\usepackage{amssymb}
\usepackage{cite}

\newtheorem{theorem}{Theorem}
\newtheorem{definition}{Definition}

\newtheorem{lemma}{Lemma}

\newtheorem{corollary}{Corollary}
\newtheorem{remark}{Remark}

\title{\LARGE \bf {Towards Finding the Critical Value for Kalman Filtering with Intermittent Observations}}

\author{Yilin Mo and Bruno Sinopoli
\thanks{
This research was supported in part by CyLab at Carnegie Mellon under grant DAAD19-02-1-0389 from the Army Research Office. Foundation. The views and conclusions contained here are those of the authors and should not be interpreted as necessarily representing the official policies or endorsements, either express or implied, of ARO, CMU, or the U.S. Government or any of its agencies.}
\thanks{
Yilin Mo and Bruno Sinopoli are with the ECE department of Carnegie Mellon University, Pittsburgh, PA {ymo@andrew.cmu.edu, brunos@ece.cmu.edu}
}
}

\begin{document}
\maketitle

\begin{abstract}
  In~\cite{Sinopoli:04}, Sinopoli et al. analyze the problem of optimal estimation for linear Gaussian systems where packets containing observations are dropped according to an i.i.d.  Bernoulli process, modeling a memoryless erasure channel. In this case the authors show that the Kalman Filter is still the optimal estimator, although boundedness of the error depends directly upon the channel arrival probability, $p$. In particular they also prove the existence of a critical value, $p_c$, for such probability, below which the Kalman filter will diverge. The authors are not able to compute the actual value of this critical probability for general linear systems, but provide upper and lower bounds. They are able to show that for special cases, i.e. $C$ invertible, such critical value coincides with the lower bound.  This paper computes the value of the critical arrival probability, under minimally restrictive conditions on the matrices $A$ and $C$. 
\end{abstract}

\section{Introduction}
\label{sec:introduction}
A large wealth of applications demand wireless communication among small embedded devices. Wireless Sensor Network (WSN) technology provides the architectural paradigm to implement systems with a high degree of temporal and spatial granularity. Applications of sensor networks are becoming ubiquitous, ranging from environmental monitoring and control to building automation, surveillance and many others\cite{wireless_sensor_network}. Given their low power nature and the requirement of long lasting deployment, communication between devices is power constrained and therefore limited in range and reliability. Changes in the environment, such as the simple relocation of a large metal object in a room or the presence of people, will inevitably affect the propagation properties of the wireless medium. Channels will be time-varying and unreliable.  Spurred by this consideration, our effort concentrates on the design and analysis of estimation and control algorithms over unreliable networks. A substantial body of literature has been devoted to such issues in the past few years. In this paper we want to revisit the paper of Sinopoli et al.~\cite{Sinopoli:04}. In that paper, the authors analyze the problem of optimal state estimation for discrete-time linear Gaussian systems, under the assumption that observations are sent to the estimator via a memoryless erasure channel. This implies the existence of a non-unitary arrival probability associated with each packet. Consequently some observations will inevitably be lost. In this case although the Kalman Filter is still the optimal estimator, the boundedness of its error depends on the arrival probabilities of the observation packets. In particular the authors prove the existence of a critical arrival probability $p_c$, below which the expectation of estimation error covariance matrix $P_k$ of Kalman filter will diverge.  The authors are not able to compute the actual value of this critical probability for general linear systems, but provide upper and lower bounds. They are able to show that for special cases such critical value coincides with the lower bound.

A significant amount of research effort has been made toward finding the critical value. In~\cite{Sinopoli:04}, the author prove that the critical value coincides with the lower bound in a special case when the system observation matrix $C$ is invertible. The condition is further weakened by Plarre and Bullo~\cite{kp-fb:07j} to $C$ only invertible on the observable subspace. In \cite{Yilin08}, the authors prove that if the eigenvalues of system $A$ matrix have distinguished absolute values, then the lower bound is indeed the critical value. The authors also provide a counter example to show that in general the lower bound is not tight. 

Other variations of the original problem are also considered. In \cite{Xu:05}, the authors introduce smart sensors, which send the local Kalman estimation instead of raw observation. In \cite{Craig:07}, a similar scenario is discussed where the sensor sends a linear combination of the current and previous measurement.  A Markovian packet dropping model is introduced in \cite{Minyi:07} and a stability criterion was given. In \cite{Liu:04}, the authors study the case where the observation at each time splits into two parts, which are sent to the Kalman filter through two independent erasure channels. A much more general model, which considered packet drop, delay and quantization of measurements in the same time, is introduced by Xie and Shi~\cite{xielihua09}.

Another interesting direction to characterize the impact of lossy network on state estimation is to directly calculate the probability distribution of estimation error covariance matrix $P_k$ instead of considering the boundedness of its expectation. In \cite{censi09}, the author gives a closed-form expression for cumulative distribution function of $P_k$ when the system satisfies non-overlapping conditions. In \cite{Ali:09}, the authors provide a numerical method to calculate the eigen-distribution of $P_k$ under the assumption that the observation matrix $C$ is random and time varying.

In the meantime, lots of research effort has been made to design estimation and control schemes over lossy network, by leveraging the result obtained from above work. In \cite{Gupta:05}, the authors consider a stochastic sensor scheduling scheme, which randomly selected one sensor to transmit observation at each time. In \cite{Ling09}, the authors shows how to design the packet arrival rate to balance the state estimation error and energy cost of packet transmission. 

In a nutshell, we feel that derivation of critical value is not only important for analyzing the performance of the system in lossy networks, but also critical for network control protocol design. However, in a large proportion of the above work, the critical value is derived under the condition that $C$ matrix is invertible or other similar conditions, which are not easy to satisfy for certain real applications\footnote{$C$ invertible implies that the number of sensors is no less than the number of states.}. In this paper, we would like to characterize the critical value under more general conditions showing that it meets the lower bound in most cases. We also study some systems for which the lower bound is not tight and try to give some insights on why this is the case.

The paper are organized in the following manner: Section~\ref{sec:formulation} formulates the problem. Section~\ref{sec:main} states all the important results of the paper, which will be proved later by Section~\ref{sec:property}, \ref{sec:critfronondeg}, \ref{sec:2degsys}. Finally Section~\ref{sec:conclusion} concludes the paper. 
\section{Problem Formulation} \label{sec:formulation}
Consider the following linear system
\begin{equation}
  \begin{split}
    x_{k + 1}  &= A x_k  + w_k , \\
    y_k  &= C x_k  + v_k  ,\\
  \end{split}
  \label{eq:systemdiscription}
\end{equation}
where $x_k\in \mathbb {R}^n$ is the state vector, $y_k \in \mathbb{R}^m$ is the output vector, $w_k \in \mathbb{R}^n$ and $v_k \in \mathbb{R}^m$ are Gaussian random vectors with zero mean and covariance matrices $Q > 0$ and $R > 0$, respectively. Assume that the initial state, $x_0$ is also a Gaussian vector of mean $\bar x_0$ and covariance matrix $\Sigma_0 > 0$. Let $w_i,v_i,x_0$ to be mutually independent. Note that we assume the covariance matrices of $w_i,v_i,x_0$ to be strictly positive definite. Define $|\lambda_1| \geq |\lambda_2| \geq \cdots \geq |\lambda_n|$ as the eigenvalues of $A$.

Consider the case where observations are sent to the estimator via a memoryless erasure channel, where their arrival is modeled by a Bernoulli independent process $\{\gamma_k\}$.  According to this model, the measurement $y_k$ sent at time $k$ reaches its destination if $\gamma_k=1$; it is lost otherwise. Let $\gamma_k$ be independent of $w_k,v_k,x_0$, i.e. the communication channel is independent of both process and measurement noises and let $P(\gamma_k=1)=p$.

The Kalman Filter equations for this system were derived in~\cite{Sinopoli:04} and take the following form:
\[
\begin{split}
  \hat x_{k|k}  &= \hat x_{k|k - 1}  + \gamma _k K_k (y_k  - C\hat  x_{k|k - 1} ) ,\\
  P_{k|k}  &= P_{k|k - 1}  - \gamma _k K_k CP_{k|k - 1}  ,\\
\end{split}
\]
where
\[
\begin{split}
  \hat x_{k + 1|k}  &= A\hat x _{k|k}  , \quad P_{k + 1|k}  = AP_{k|k} A^T  + Q ,\\
  K_k  &= P_{k|k - 1} C^T (CP_{k|k - 1} C^T  + R)^{ - 1}  ,\\
  \hat x_{0| - 1}  &= \bar x_0 ,\quad P_{0| - 1}  = \Sigma _0  .\\
\end{split}
\]
In the hope to improve the legibility of the paper we will slightly abuse the notation, by substituting $P_{k|k-1}$ with $P_k$. The equation for the error covariance of the one-step predictor is the following:
\begin {equation}
\label{eq:recursive}
P_{k+1}=AP_k A^T+Q-\gamma_k A P_k C^T(C P_k C^T+R)^{-1}C P_k A^T.
\end {equation}

If $\gamma_k$s are i.i.d. Bernoulli random variables, the following theorem holds~\cite{Sinopoli:04}:
\begin{theorem}\label{theorem:lambdaCrit}
  If $(A,Q^{\frac{1}{2}})$ is controllable, $(C,A)$ is detectable, and $A$ is unstable, then there exists a $p_c\in[0,1)$ such that \footnote{We use the notation $\sup_{k}A_k=+\infty$ when the sequence $A_k\geq 0$ is not bounded; i.e., there is no matrix $M\geq 0$ such that $A_k\leq M, \forall k$.}$^,$\footnote{Note that all the comparisons between matrices in this paper are in the sense of positive definite if without further notice}
  \begin{eqnarray}\label{eqn:lambdaCrit}
    \sup_{k}EP_k=+\infty & \mbox{for } 0\leq p \leq p_c \; \; \mbox{and  \hspace{3mm}$\exists P_0 \geq 0$},&\\   
    EP_k \leq M_{P_0} \;\; \forall t & \mbox{for } p_c < p\leq 1 \; \; \mbox{and \hspace{3mm}$\forall P_0 \geq 0$}, &
  \end{eqnarray}
  where $M_{P_0} > 0 $  depends on the initial condition $P_0 \geq 0$.
\end{theorem}

For simplicity, we will say that $EP_k$ is unbounded if $\sup_{k} EP_k = +\infty$ or $EP_k$ is bounded if there exists a uniform bound independent of $k$.

\section{Main Result}\label{sec:main}

In this section, we want to state all the important results for critical value, the proof of which can be found in later sections. Through out the rest of the paper, we always assume that the following conditions hold:
\begin{enumerate}
  \item[$(H1)$] $(C,A)$ is detectable.
  \item[$(H2)$] $A$ can be diagonalized.
  \item[$(H3)$] $R , Q , \Sigma_0$ are strict positive definite.
\end{enumerate}

From Section~\ref{sec:formulation}, it is clear that the critical value of a system should be a function of all system parameters, i.e. $A,\,C,\,Q,\,R,\,\Sigma_0$. However, the following theorem, the proof of which is in Section~\ref{sec:property}, states that the critical value does not depend on $R,\,Q,\,\Sigma_0$ as long as they are all positive definite. 

\begin{theorem}
  \label{theorem:independent}
  If $R,Q,\Sigma_0 > 0$ are strictly positive definite, then the critical value of a system is just a function of $A,C$, and is independent of $R , Q , \Sigma_0$.
\end{theorem}
Since we have already assumed that $R,Q,\Sigma_0$ are strictly positive definite, by Theorem~\ref{theorem:independent}, we can let $R = I_m, Q = I_n, \Sigma_0 = I_n$ without loss of generality. Also since we assume that $A$ can be diagonalized, we can always transform the system into its diagonal standard form. Hence, we assume that $A$ is diagonal. We can also denote $f_c(A,C)$ as the critical value of system $(A,C)$. 

When the dimension of $A,\,C$ is large, which is often the case in reality, it is desirable to break the large system to several smaller blocks (or subsystems), which are easier to analyse. As a result, we define a block of the system in the following way:
\begin{definition}
  Consider the system $(A,C)$ is in its diagonal standard form, which means $A = diag(\lambda_{1},\ldots,\lambda_{n})$ and $C = [C_1,\ldots,C_n]$. A block of the system is defined as subsystem $(A_{\mathcal I} = diag(\lambda_{i_1},\ldots,\lambda_{i_l}),C_{\mathcal I} = [C_{i_1},\ldots,C_{i_l}]),\;1\leq i_1 < \ldots < i_l \leq n$, where $\mathcal I = \{i_1,\ldots,i_l\}\subset \{1,\ldots,n\}$ is the index set.
\end{definition}
A special type of block, which we call equi-block, plays a central role in determining the critical value of the system and it is defined as
\begin{definition}
  An equi-block is a block which satisfies $|\lambda_{i_1}| = \ldots = |\lambda_{i_l}|$, and we denote it as $(A_{\mathcal I_e},C_{\mathcal I_e})$, where $\mathcal I_e$ is the index set. 
\end{definition}
\begin{definition}
  $\mathcal D(A,C)$ is defined as the dimension of the largest equi-block of the system. 
\end{definition}
The following theorem shows a basic inequality between the critical value of the original system and smaller blocks, which we will prove in Section~\ref{sec:property}.
\begin{theorem}
  \label{theorem:criticalfunction}
  Define $f_c(A,C)$ as the critical value for system $(A,C)$. If $A = diag(\lambda_1,\ldots,\lambda_n)$ diagonal and $C = [C_1,\ldots,C_n]$, then 
  \begin{equation}
    f_c(A,C) \geq  f_c(A_{\mathcal I},C_{\mathcal I}),
  \end{equation}
  for all possible index set $\mathcal I\subset \{1,\ldots,n\}$.
\end{theorem}

Before we continue on, we need to define the following terms:
\begin{definition}
  A system $(A,C)$ is one step observable if $C$ is full column rank.
\end{definition}

\begin{definition}
  An equi-block is degenerate if it is not one step observable. It is non-degenerate otherwise.
\end{definition}

\begin{definition}
  The system is non-degenerate if every equi-block of the system is non-degenerate. It is degenerate if there exists at least one degenerate equi-block. 
\end{definition}

For example, if $A = diag (2,-2)$ and $C =[1,1]$, then the system is degenerate since it is an equi-block and not one step observable. For $A = diag(2,-2,3,-3)$ and $C = \left[ {\begin{array}{*{20}c}
  1 & 0 & 1 & 0  \\
  0 & 1 & 0 & 1  \\
\end{array}} \right]$, two equi-blocks are $(diag(2,-2),I_2)$ and $(diag(3,-3),I_2)$ and both of them are one step observable. Thus, the system is non-degenerate. 

It can be seen that non-degeneracy is a stronger property than observability but much more weaker than one step observability. In fact, for a one step observable system, $C$ matrix must have at least $n$ rows, which implies $y_k$ is at least a $\mathbb R^n$ vector. On the other hand, for non-degenerate system, the $C$ matrix can only have $\mathcal D(A,C)$ rows. In reality, $\mathcal D(A,C)$ is usually a small number comparing to $n$. 

In \cite{Sinopoli:04}, the authors proved that the critical value meets the lower bound when the system is one step observable.  In this paper, we weaken the condition from one step observability to non-degeneracy.
\begin{theorem}
  \label{theorem:critnondeg}
  If the system \ref{eq:systemdiscription} satisfies assumptions $(H1)-(H3)$ and the equiblocks of $A$ associated to the unstable and critically stable eigenvalues are non-degenerate, the critical value of the Kalman filter is
  \begin{equation}
    p_c = \max(1- |\lambda_1|^{-2}, 0)
  \end{equation}
  where $\lambda_1$ is the dominant eigenvalue.
\end{theorem}

For degenerate systems we can show that in general the critical value is larger than the one computed in theorem \ref{theorem:critnondeg}. Nonetheless in this paper we will compute the critical value for second order degenerate systems. This includes a very practical case, involving complex conjugate eigenvalues. Let $A=diag(\lambda_1,\lambda_2)\in \mathbb R^{2\times 2}$. We can use the following theorem in conjunction with Theorem~\ref{theorem:criticalfunction} as the building block to allow analysis of larger systems.

\begin{theorem}
  \label{theorem:critdeg}
  For a detectable system with $A = diag(\lambda_1,\lambda_2),\,|\lambda_1| \geq |\lambda_2|$ and $R,\,Q,\,\Sigma_0 > 0$, the critical value is
  \begin{equation}
    p_c = f_c(A,C) = \max(1-|\lambda_1|^{-2},0),
  \end{equation}
  if the system is non-degenerate, or in other word, if one of the following conditions holds
  \begin{enumerate}
    \item $|\lambda_1| > |\lambda_2|$,
    \item $rank(C) = 2$.
  \end{enumerate}
  Otherwise the system is degenerate and its critical value is
  \begin{equation}
    p_c = f_c(A,C) =\max(1-|\lambda_1|^{-\frac{2}{1-D_M(\varphi/2\pi)}},0),
  \end{equation}
  where $\lambda_1 = \lambda_2 \exp(j\varphi)$, and $D_M(x)$ is the modified Dirichlet function defined as
  \begin{equation}
    D_M(x) = \left\{ \begin{array}{*{20}l}
      0 & \textrm{for $x$ irrational}\\
      1/q & \textrm{for $x = r/q$, $r,q \in \mathbb Z$ and irreducible.}\\
    \end{array}\right..
  \end{equation}
\end{theorem}

\section{Properties of Critical Value}
\label{sec:property}
In this section, we will prove Theorem~\ref{theorem:independent} and \ref{theorem:criticalfunction}, which demonstrate the relationship between critical value and system parameters. Throughout this section, we always assume that assumption $(H1)-(H3)$ holds.

First we want to prove the independence between critical value and the covariance matrix of the noise. 
\begin{IEEEproof}[Proof of Theorem~\ref{theorem:independent}]
  Since $R,\Sigma_0,Q >0$, we can find uniform upper and lower bounds $\underline \alpha,\overline \alpha > 0$, such that 
  \[
  \underline \alpha I_m \leq R \leq \overline \alpha I_m, \, \, \, \, \underline \alpha I_n \leq \Sigma_0 \leq \overline \alpha I_n,\, \, \, \,\underline \alpha I_n \leq Q \leq \overline \alpha I_n.\] 

  Let us define $\underline P_0 = \underline \alpha I_n,\,\overline P_0 = \overline \alpha I_n,\,P_0^* = I_n$, and 
  \begin{align*}
 \underline P_{k+1}&=A\underline P_k A^T+\underline \alpha I_n-\gamma_k A \underline P_k C^T(C \underline P_k C^T+\underline \alpha I_m)^{-1}C \underline P_k A^T,\\
 \overline P_{k+1}&=A\overline P_k A^T+\overline \alpha I_n-\gamma_k A \overline P_k C^T(C \overline P_k C^T+\overline \alpha I_m)^{-1}C \overline P_k A^T,\\
 P_{k+1}^*&=AP_k^* A^T+ I_n-\gamma_k A P_k^* C^T(C P_k^* C^T+ I_m)^{-1}C P_k^* A^T.
  \end{align*}
  By induction, it is easy to check that $\underline P_k = \underline \alpha P_k^*$ and $\overline P_k = \overline \alpha P_k^*$ for all $k$. 

  Also we know that $\underline P_0 \leq P_0$. By induction, suppose that $\underline P_k \leq P_k$, then
  \begin{align*}
 \underline P_{k+1}&=A\underline P_k A^T+\underline \alpha I_n-\gamma_k A \underline P_k C^T(C \underline P_k C^T+\underline \alpha I_m)^{-1}C \underline P_k A^T \\
 &\leq A P_k A^T+\underline \alpha I_n-\gamma_k A  P_k C^T(C  P_k C^T+\underline \alpha I_m)^{-1}C  P_k A^T \\
 &\leq A P_k A^T+Q-\gamma_k A  P_k C^T(C  P_k C^T+R)^{-1}C  P_k A^T = P_{k+1}.
  \end{align*}
Hence, $\underline P_k \leq P_k$ for all $k$. By the same argument, $P_k \leq \overline P_k$ for all $k$, which implies that
\begin{equation}
    \underline P_k = \underline \alpha P_k^* \leq P_k \leq \overline P_k=\overline \alpha P_k^*. 
\end{equation}
Since $\underline \alpha,\,\overline \alpha > 0$, the boundedness of $P_k$ is equivalent to the boundedness of $P_k^*$. However by the definition of $P_k^*$, we know that it is only a function of $A,\,C$ which is independent of $R,\,Q,\,\Sigma_0$.
\end{IEEEproof}

We now want to prove that the critical value of a system is larger than the critical value of any of its blocks.
\begin{IEEEproof}[Proof of Theorem~\ref{theorem:criticalfunction}]
	With out loss of generality\footnote{If $\mathcal I = \{1,\ldots,m\} $, the proof is trivial. If $\mathcal I$ is an arbitrary subset of size $l$, we can always permute the states to make it equal $\{1,\ldots,l\}$}, we assume that $\mathcal I = \{1,\ldots,l\}$, $l < m$. Let us define $\mathcal J = \{l+1,\ldots,m\}$ to be the complement index set of $\mathcal I$. 

  By Theorem~\ref{theorem:independent}, we suppose for the original system $R = I_m,\,Q = \Sigma_0 = I_n$.

  Let us define $\tilde P_0 = \Sigma_0 = I_n$ and 
  \begin{displaymath}
    \tilde P_{k+1}=A\tilde P_k A^T+ I_n-\gamma_k A \tilde P_k \tilde C^T(\tilde C \tilde P_k \tilde C^T+ I_{2m}/2)^{-1}\tilde C \tilde P_k A^T,
  \end{displaymath}
  where $\tilde C = \left[ {\begin{array}{*{20}c}
    C_{\mathcal I} & 0   \\
    0 & C_{\mathcal J}   \\
  \end{array}} \right]\in \mathbb R^{2m\times n}$. Using Matrix Inversion Lemma, we can show that 
  \begin{align}
    \label{eq:informationfilter}
    P_{k+1} &=A\left(P_k^{-1} + \gamma_k C^T C \right)^{-1}A^T + I_n .\\
    \tilde P_{k+1} &=A\left[\tilde P_k^{-1} + 2\gamma_k \tilde C^T \tilde C \right]^{-1}A^T + I_n .\nonumber
  \end{align}
  We know that $P_0 =\tilde  P_0$. By induction, suppose that $ P_k \geq \tilde P_k$, then 
  \begin{equation}
    \label{eq:subsysteminequality}
    \begin{split}
     & P_k^{-1} + \gamma_k C^TC -\left( \tilde P_k^{-1} + 2\gamma_k\tilde C^T \tilde C\right) = P_{k}^{-1}-\tilde P_k^{-1} + \gamma_k\left(\left[ {\begin{array}{*{20}c}
	C_1^TC_1 & C_1^T C_2\\
	C_2^TC_1 & C_2^TC_2\\
      \end{array}}\right] - 2\left[ {\begin{array}{*{20}c}
	C_1^TC_1 & 0\\
	0 & C_2^TC_2\\
      \end{array}}\right]\right)\\
      &= P_{k}^{-1}-\tilde P_k^{-1}+\gamma_k\left[ {\begin{array}{*{20}c}
	- C_1^TC_1 & C_1^T C_2\\
	C_2^TC_1 &  - C_2^TC_2\\
      \end{array}}\right] =P_{k}^{-1}-\tilde P_k^{-1}-\gamma_k\left[ {\begin{array}{*{20}c}
	C_1 & -C_2\\
      \end{array}}\right]\left[ {\begin{array}{*{20}c}
	C_1^T\\
	-C_2^T\\
      \end{array}}\right] \leq 0 .
    \end{split}
  \end{equation}
  Hence,
  \[
  P_{k+1}  = A\left(P_k^{-1} + \gamma_k C^T  C \right)^{-1}A^T + Q  \geq A\left(\tilde P_k^{-1} + 2\gamma_k\tilde  C^T  \tilde C \right)^{-1}A^T + Q = \tilde P_{k+1}.
  \]
  Thus, by induction, $P_k \geq \tilde P_k, \forall k$, which in turn proves 
  \begin{equation}
  f_c(A,C) \geq f_c(A,\tilde C).
    \label{eq:criticalfun1}
  \end{equation}

  Now define $\tilde P_{0,\mathcal I} = I_l$, $\tilde P_{0,\mathcal J} = I_{m-l}$ and
  \begin{align*}
	  \tilde P_{k+1,\mathcal I}&=A_{\mathcal I}\tilde P_{k,\mathcal I} A_{\mathcal I}^T+ I_l-\gamma_k A_{\mathcal I} \tilde P_{k,\mathcal I}  C_{\mathcal I}^T(C_{\mathcal I} \tilde P_{k,\mathcal I} C_{\mathcal I}^T+ I_m/2)^{-1}C_{\mathcal I} \tilde P_{k,\mathcal I} A_{\mathcal I}^T,\\
	  \tilde P_{k+1,\mathcal J}&=A_{\mathcal J}\tilde P_{k,\mathcal J} A_{\mathcal J}^T+ I_{m-l}-\gamma_k A_{\mathcal J} \tilde P_{k,\mathcal J}  C_{\mathcal J}^T(C_{\mathcal J} \tilde P_{k,\mathcal J} C_{\mathcal J}^T+ I_m/2)^{-1}C_{\mathcal J} \tilde P_{k,\mathcal J} A_{\mathcal J}^T.\\
  \end{align*}
  It is not hard to check that $\tilde P_k = \left[ {\begin{array}{*{20}c}
    \tilde P_{k,\mathcal I} & 0   \\
    0 & \tilde P_{k,\mathcal J}   \\
  \end{array}} \right]$, for all $k$. As a result, $\tilde P_k$ is bounded if and only if $\tilde P_{k,\mathcal I}$ and $\tilde P_{k,\mathcal J}$ are both bounded. Combining \eqref{eq:criticalfun1}, we know
  \begin{displaymath}
    f_c(A,C) \geq f_c(A,\tilde C) = \max \{f_c(A_{\mathcal I},C_{\mathcal I}),f_c(A_{\mathcal J},C_{\mathcal J})\}.
  \end{displaymath}
\end{IEEEproof}

\section{Critical Value for Non-degenerate Systems}
\label{sec:critfronondeg}
This section is devoted to proving Theorem~\ref{theorem:critnondeg}. Before continuing, we would like to state several important intermediate results which are useful for proving the main theorem and have some theoretical value of their own.

We will first deal with systems whose eigenvalues are all unstable. We will lift this restriction later in the paper. By ``unstable'' we mean that its absolute value is strictly greater than $1$. We will call the eigenvalues on the unit circle critically stable and the ones with absolute values strictly less than $1$ stable. Since $A$ is diagonalizable, we will restrict our analysis to systems with diagonal $A$. Also since some eigenvalues of $A$ may be complex, we will use Hermitian instead of transpose in the remaining part of the article. 

Similarly to the observability Grammian, we find that the matrix $\sum_{i=1}^\infty \gamma_i A^{-iH}C^HCA^{-i}$ plays an essential role in determining the boundedness of Kalman filter, which is characterized by the following theorem:

\begin{theorem}
  \label{theorem:uniformbound}
  If a system satisfies assumptions $(H1)-(H3)$ and all its eigenvalues are unstable, then $\sup_k EP_k$ if finite if and only if $E[(\sum_{i=1}^\infty \gamma_i A^{-iH}C^HCA^{-i})^{-1}]$ exists, and it satisfies the following inequality
  \begin{equation}
    \underline \alpha E[(\sum_{i=1}^\infty \gamma_i A^{-iH}C^HCA^{-i})^{-1}] \leq \sup_k EP_k \leq \overline \alpha E[(\sum_{i=1}^\infty \gamma_i A^{-iH}C^HCA^{-i})^{-1}],
  \end{equation}
  where $\underline \alpha, \overline \alpha > 0$ are constants.
\end{theorem}

By manipulating  $\sum_{i=1}^\infty \gamma_i A^{-iH}C^HCA^{-i}$, we established the following result, which is essentially equivalent to Theorem~\ref{theorem:critnondeg}, but restricted on the systems whose eigenvalues are all unstable.

\begin{theorem}
  \label{theorem:bound} 
  If a non-degenerate system satisfies assumptions $(H1)-(H3)$ and all its eigenvalues are unstable, then the critical value of the system is
  \[
  p_c=1-|\lambda_1|^{-2}.
  \]
  If the arrival probability $p \geq p_c$, then for all initial conditions, $EP_k$ will be bounded for all $k$. Else if $p\leq p_c$, for some initial conditions, $EP_k$ is unbounded.
\end{theorem}

Now we need to generalize this result to systems that have stable eigenvalues. The following theorem provides an important inequality for systems that have stable eigenvalues:

\begin{theorem}
  \label{theorem:relaxation}
  Consider a system satisfies assumption $(H1) - (H3)$ with a diagonal $A = diag (A_1, A_2 ,A_3)$, $C = [C_1,C_2,C_3]$\footnote{Note that there is no requirement of non-degeneracy for this theorem.}. If $A_1$ is the unstable part, $A_2$ is the critically stable part and $A_3$ is the stable part and $C_i$s are of proper dimensions, then the critical value of the system satisfies the following inequality
  \begin{equation}
	  f_c(A,C) \leq \lim_{\alpha \rightarrow 1^+} f_c(diag(\alpha A_1,\alpha A_2),[C_1,C_2])\footnote{Note that $diag(\alpha A_1,\alpha A_2)$ are unstable when $\alpha > 1$. Hence the right hand side of the inequality can be computed by Theorem~\ref{theorem:bound}.}.
    \label{eq:stabledoesnotmatter}
  \end{equation}
\end{theorem}

Combining Theorem~\ref{theorem:criticalfunction}, \ref{theorem:bound} and \ref{theorem:relaxation}, we will prove Theorem~\ref{theorem:critnondeg} in the last part of this section.

\subsection{Proof of Theorem~\ref{theorem:uniformbound}}

In this subsection, we are going to prove Theorem~\ref{theorem:uniformbound}. The key idea in the proof is to avoid analysing Riccati Equations, which were intensively studied in the previous works\cite{Sinopoli:04}\cite{kp-fb:07j}, and try to formulate estimation error covariance $P_k$ by maximum likelihood estimator. First let us write down the relation between $\gamma_i y_i$ and $x_k$: 
\begin{equation}
  \label{eq:relationloss}
  \left[ {\begin{array}{*{20}c}
    {\gamma_k y_k } \\
    \vdots \\
    {\gamma_0y_0 } \\
    \bar x_0 \\
  \end{array}} \right] = \left[ {\begin{array}{*{20}c}
    {\gamma_k CA^{ - 1} } \\
    \vdots \\
    {\gamma_{0}CA^{ - k-1} } \\
    {A^{ - k - 1} } \\
  \end{array}} \right]x_{k + 1} + \left[ {\begin{array}{*{20}c}
    {\gamma_{k}v_k } \\
    \vdots \\
    {\gamma_{0}v_0 } \\
    { \bar x_0 - x_0 } \\
  \end{array}} \right] - \left[ {\begin{array}{*{20}c}
    {\gamma_k CA^{ - 1} } & \cdots & 0 & 0 \\
    \vdots & \ddots & \vdots & \vdots \\
    {\gamma_{0}CA^{ - k-1} } & \cdots & {\gamma_{0}CA^{ - 2} } & \gamma_0 CA^{-1} \\
    {A^{ - k - 1} } & \cdots & {A^{ - 2} } & {A^{ - 1} } \\
  \end{array}} \right] \left[ {\begin{array}{*{20}c}
    {w_k } \\
    \vdots \\
    {w_1 } \\
    {w_0 } \\
  \end{array}} \right].
\end{equation}
The rows where $\gamma_i$s are zero can be deleted, since they do not provide any information to improve the estimation of $x_{k+1}$. To write \eqref{eq:relationloss} in a more compact way, let us define the following quantities:
\begin{equation}
  F_k \triangleq\left[ {\begin{array}{*{20}c}
    {A^{ - 1} } & \cdots & 0 & 0 \\
    \vdots & \ddots & \vdots & \vdots \\
    {A^{ - k} } & \cdots & {A^{ - 1} } & 0 \\
    {A^{ - k-1} } & \cdots & {A^{ - 2} } & {A^{ - 1} } \\
  \end{array}} \right] \in \mathbb{R}^{n(k+1)\times n(k+1)}.\\
  \label{eq:matrixdefF}
\end{equation}
\begin{equation}
  G_k\triangleq \in \left[ {\begin{array}{*{20}c}
    {C} & \cdots & 0 \\
    \vdots & \ddots & \vdots \\
    {0 } & \cdots & C \\
    {0 } & \cdots & I_n  \\
  \end{array}} \right] \mathbb{R}^{[n+m(k+1)]\times n(k+1)}.
  \label{eq:matrixdefG}
\end{equation}
\begin{equation}
  e_k\triangleq -G_k F_k\left[ {\begin{array}{*{20}c}
    {w_k } \\
    \vdots \\
    {w_1 } \\
    {w_0 } \\
  \end{array}} \right] + \left[ {\begin{array}{*{20}c}
    {v_k } \\
    \vdots \\
    {v_0 } \\
    { \bar x_0- x_0 } \\
  \end{array}} \right] \in \mathbb{R}^{n+m(k+1)}.
  \label{eq:matrixdefe}
\end{equation}
\begin{equation}
  T_k \triangleq \left[ {\begin{array}{*{20}c}
    {CA^{ - 1} } \\
    \vdots \\
    {CA^{ - k-1} } \\
    {A^{ - k - 1} } \\
  \end{array}} \right] \in \mathbb{R}^{(n+m(k+1))\times n},
  Y_k = \left[ {\begin{array}{*{20}c}
    {y_k } \\
    \vdots \\
    {y_0 } \\
    \bar x_0 \\
  \end{array}} \right] \in \mathbb{R}^{n+m(k+1)}.
\end{equation}

Define $\Gamma_k$ as the matrix of all non zero rows of $diag(\gamma_kI_m,\gamma_{k-1}I_m,\ldots,\gamma_0I_m,I_n)$.  Thus $\Gamma_k$ is a ($m\sum_{i=0}^k \gamma_i+n$) by ($n+m(k+1)$) matrix. Also define

\[
\widetilde Y_k\triangleq \Gamma_k Y_k,\,\widetilde T_k\triangleq \Gamma_kT_k,\,\widetilde e_k\triangleq \Gamma_ke_k.
\]
$\widetilde Y_k,\,\widetilde T_k,\,\widetilde e_k$ are now stochastic matrices as they are functions of $\gamma_k,\gamma_{k-1},\ldots,\gamma_0$.

We can rewrite \eqref{eq:relationloss} in a more compact form as
\begin{equation}
  \label{eq:linearform}
  \widetilde Y_k=\widetilde T_k x_{k+1}+\widetilde e_k.
\end{equation}

From \eqref{eq:linearform}, we know that $\widetilde Y_k$ is Gaussian distributed with unknown mean $\widetilde T_k x_{k+1}$ and known covariance $Cov(\widetilde e_k|\Gamma_k)$. Hence, we can prove the following lemmas:

\begin{lemma}
  \label{lemma:kalmanerror}
  If $A$ is invertible, then the state estimation $\hat x_{k+1|k}$ and estimation error covariance $P_{k+1}$ given by Kalman filter satisfy the following equations
  \begin{equation}
    \hat x_{k+1|k}=(\widetilde T_k^H Cov (\widetilde e_k|\Gamma_k)^{-1}\widetilde T_k)^{-1}T_k^H Cov(\widetilde e_k|\Gamma_k)^{-1}\widetilde Y_k,
    \label{eq:kalmangain}
  \end{equation}
  \begin{equation}
    P_{k+1}=(\widetilde T_k^HCov(\widetilde e_k|\Gamma_k)^{-1}\widetilde T_k)^{-1}.
    \label{eq:kalmanerror}
  \end{equation}
\end{lemma}

\begin{lemma}
  \label{lemma:bound}
  If $A = diag (\lambda_1,\lambda_2,\cdots,\lambda_n)$, where $|\lambda_1| \geq |\lambda_2| \geq \cdots \geq |\lambda_n| > 1$, then $F_k F_k^H$ is bounded by
  \begin {equation}
  \frac{1}{(|\lambda_1|+1)^2} I_{n(k+1)} \leq F_k F_k^H \leq \frac{1}{(|\lambda_n|-1)^2} I_{n(k+1)},
\end {equation}
where $F_k$ is defined in \eqref{eq:matrixdefF}.
\end{lemma}

\begin{lemma}
  \label{lemma:pseudoinv}
  If a system satisfies assumptions $(H1)-(H3)$ and all its eigenvalues are unstable, then the error covariance matrix of Kalman Filter is bounded by
  \begin {equation}
  \label{eq:mlebound}
  \underline \alpha (\widetilde T_k^H\widetilde T_k)^{-1} \leq P_{k+1} \leq \overline \alpha(\widetilde T_k^H\widetilde T_k)^{-1}
\end {equation}
where $\underline \alpha,\overline \alpha \in \mathbb {R}$ are constants independent of $\gamma_i$ and $k$.\footnote{We abuse the notations of $\overline \alpha,\underline \alpha$, which will also be used several times later. These notations only means constant lower and upper bound, which are not necessarily the same in different theorems.}
\end{lemma}
\begin{IEEEproof}[Proof of Lemma~\ref{lemma:kalmanerror}]
  Given the observation $\widetilde Y_k$, by \eqref{eq:linearform}, we know that the maximum likelihood estimator of $x_{k+1}$ is 
  \begin{displaymath}
    \hat x_{k+1|k}=(\widetilde T_k^H Cov (\widetilde e_k|\Gamma_k)^{-1}\widetilde T_k)^{-1}T_k^H Cov(\widetilde e_k|\Gamma_k)^{-1}\widetilde Y_k,
  \end{displaymath}
  and the estimation error covariance is
  \begin{displaymath}
    P_{k+1}=(\widetilde T_k^HCov(\widetilde e_k|\Gamma_k)^{-1}\widetilde T_k)^{-1}.
  \end{displaymath}
  Since $Y_k$ is Gaussian with unknown mean $\widetilde T_k x_{k+1}$ and known covariance $Cov(\widetilde e_k|\Gamma_k)$, we know that the maximum likelihood estimator is the optimal one in minimum mean error covariance sense. Thus, $\hat x_{k+1|k}$ and $P_k$ given by maximum likelihood estimator are essentially the same as $\hat x_{k+1|k}$ and $P_k$ of Kalman filter, which concludes the proof. 
\end{IEEEproof}

\begin{IEEEproof}[Proof of Lemma~\ref{lemma:bound}]
  Notice that
  \[
  F_k^{-1}=
  \left[ {\begin{array}{*{20}c}
    A & {} & {} & {} \\
    { - I} & \ddots & {} & {} \\
    {} & \ddots & A & {} \\
    {} & {} & { - I} & A \\
  \end{array}} \right].
  \]
  Therefore,
  \[
  (F_k F_k^H)^{-1}= \left[ {\begin{array}{*{20}c}
    {AA^H + I} & { - A} & {} & {} \\
    { - A^H} & \ddots & \ddots & {} \\
    {} & \ddots & {AA^H + I} & { - A} \\
    {} & {} & { - A^H} & {AA^H } \\
  \end{array}} \right].
  \]
  By Gershgorin's Circle Theorem\cite{gerschgorincircle}, we know that all the eigenvalues of $(F_kF_k^H)^{-1}$ are located inside one of the following circles:$|\zeta - |\lambda_i|^2-1| = |\lambda_i|,|\zeta - |\lambda_i|^2-1| = 2|\lambda_i|,|\zeta - |\lambda_i|^2| = |\lambda_i|$, where $\zeta$s are the eigenvalues of $(F_k F_k^H)^{-1}$.

  Since $|\lambda_1| \geq |\lambda_2| \geq \cdots \geq |\lambda_n| \geq 1$, for each eigenvalue of $(F_kF_k^H)^{-1}$, the following holds:
  \begin{equation}
    \zeta \geq min \{|\lambda_i|^2+1-|\lambda_i|,|\lambda_i|^2+1-2|\lambda_i|,|\lambda_i|^2-|\lambda_i|\},
  \end {equation}
  and
  \begin{equation}
    \zeta \leq \max \{|\lambda_i|^2+1+|\lambda_i|,|\lambda_i|^2+1+2|\lambda_i|,|\lambda_i|^2+|\lambda_i|\}.
  \end {equation}

  Thus, $0 < (|\lambda_n|-1)^2 \leq \zeta \leq (|\lambda_1|+1)^2$, which in turn gives
  \[
  \frac{1}{(|\lambda_1|+1)^2}I_{n(k+1)} \leq F_k F_k^H \leq \frac{1}{(|\lambda_n|-1)^2}I_{n(k+1)}.
  \]
\end{IEEEproof}

\begin{IEEEproof}[Proof of Lemma~\ref{lemma:pseudoinv}]
  Since $w_i,v_j,x_0$ are mutually independent,
  \[
  \begin {split}
  Cov (e_k) & = Cov(G_k F_k[w_k,\ldots,w_0]^T)+Cov([v_k,\ldots,v_1,x_0]) \\
  & =G_k F_k diag(Q,Q,\ldots,Q) F_k^H G_k^H +diag(R,R,\ldots,R,\Sigma_0).
\end{split}
\]

Since we assume that $Q = I_m, R = \Sigma_0 = I_n$, using Lemma~\ref{lemma:bound}, it is easy to show that
\[
\frac{1}{(|\lambda_1|+1)^2}G_k G_k^H+I_{n+mk} \leq Cov(e_k) = G_k F_k F_k^HG_k^H+I_{n+mk} \leq \frac{1}{(|\lambda_n|-1)^2}G_k G_k^H+I_{n+mk}.
\]
Since $G_k G_k^H = diag(CC^H,\ldots,CC^H,CC^H + I_n)$, define $\underline n_G = \lambda_{min}(CC^H)$ and $\overline n_G = \lambda_{max}(CC^H) + 1$ , we know that $\underline n_G I_{n+mk} \leq G_k G_k^H \leq \overline n_G I_{n+mk}$, which implies
\[
\underline \alpha \Gamma_k\Gamma_k^T \leq Cov(\widetilde e_k|\Gamma_k) =\Gamma_kCov(e_k)\Gamma_k^T\leq \overline \alpha \Gamma_k\Gamma_k^T,
\]
where $\underline \alpha = \frac{\underline n_G}{(|\lambda_1|+1)^2}+1,\overline \alpha =\frac{\overline n_G}{(|\lambda_n|-1)^2}+1$.  Notice that $\Gamma_k\Gamma_k^T=I$. Therefore,
\[
\underline \alpha I \leq Cov(\widetilde e_k|\Gamma_k) \leq \overline \alpha I.
\]
The above bound is independent of $k$ and $\gamma_i$, which proves
\[
\underline \alpha (\widetilde T_k^H\widetilde T_k)^{-1} \leq P_{k+1}
= (\widetilde T_k^HCov(\widetilde e_k|\Gamma_k)^{-1}\widetilde T_k)^{-1} \leq
\overline \alpha(\widetilde T_k^H\widetilde T_k)^{-1}.
\]
\end{IEEEproof}

Now we are ready to prove Theorem~\ref{theorem:uniformbound}.
\begin{IEEEproof}[Proof of Theorem~\ref{theorem:uniformbound}]
As a result of Lemma~\ref{lemma:pseudoinv}, we only need to show that
\begin{equation}
	\underline \alpha E\left[(\sum_{i=1}^\infty \gamma_i A^{-iH}C^HCA^{-i})^{-1}\right] \leq \sup_k E\left[\left(\widetilde T_k^H \widetilde T_k\right)^{-1}\right] \leq \overline \alpha E\left[(\sum_{i=1}^\infty \gamma_i A^{-iH}C^HCA^{-i})^{-1}\right]
\end{equation}
Rewrite $\widetilde T_k^H\widetilde T_k$ as
  \begin{equation}
    \begin{split}
      \widetilde T_k^H\widetilde T_k = \sum_{i=1}^{k+1}\gamma_{k+1-i} A^{-iH}C^HCA^{-i}+A^{-(k+1)H}A^{-(k+1)}.
    \end{split}
  \end{equation}
  We know that $\widetilde T_k^H\widetilde T_k$ will have the same distribution as 
  \[
  \sum_{i=1}^{k+1}\gamma_i A^{-iH}C^HCA^{-i}+A^{-(k+1)H}A^{-(k+1)},
  \]
  since $\gamma_i$s are i.i.d. distributed.

  First we want to prove the right hand side of the inequality. Since all the eigenvalues of $A$ are unstable, 
  \begin{displaymath}
    \sum_{i=1}^{\infty} A^{-iH}C^HCA^{-i} \leq \beta I,
  \end{displaymath}
  where $\beta>0 $ is a constant. Thus
  \begin{displaymath}
    \begin{split}
      &\sum_{i=1}^{k+1}\gamma_i A^{-iH}C^HCA^{-i}+A^{-(k+1)H}A^{-(k+1)} \\
      &\geq \sum_{i=1}^{k+1}\gamma_i A^{-iH}C^HCA^{-i}+A^{-(k+1)H} \left(\beta^{-1} \sum_{i=1}^{\infty} A^{-iH}C^HCA^{-i}\right) A^{-(k+1)}\\
      &=\sum_{i=1}^{k+1}\gamma_i A^{-iH}C^HCA^{-i}+\beta^{-1} \sum_{i=k+2}^{\infty} A^{-iH}C^HCA^{-i}\\
      &\geq\sum_{i=1}^{k+1}\gamma_i A^{-iH}C^HCA^{-i}+\beta^{-1} \sum_{i=k+2}^{\infty}\gamma_i A^{-iH}C^HCA^{-i}\geq \min (1,\beta^{-1}) \sum_{i=1}^{\infty} \gamma_i A^{-iH}C^HCA^{-i}.
    \end{split}
  \end{displaymath}
  Thus,
  \begin{equation}
    E \left[(\widetilde T_k^H \widetilde T_k)^{-1}\right] \leq  \max (1,\beta) E\left[(\sum_{i=1}^{\infty} \gamma_i A^{-iH}C^HCA^{-i})^{-1}\right], \quad \textrm{for all $k$},
  \end{equation}
  which proves the right hand side of the inequality.

  Then we want proof the left hand side of the inequality. We know that
  \begin{displaymath}
    \sum_{i=1}^{k}\gamma_i A^{-iH}C^HCA^{-i}+A^{-(k+1)H}A^{-(k+1)} \leq \sum_{i=1}^{\infty}\gamma_i A^{-iH}C^HCA^{-i}+A^{-(k+1)H}A^{-(k+1)}.
  \end{displaymath}
  Thus, 
  \begin{displaymath}
    \begin{split}
      \sup_k E[(\widetilde T_k^H \widetilde T_k)^{-1}]& = \sup_{k} E[(\sum_{i=1}^{k+1}\gamma_i A^{-iH}C^HCA^{-i}+A^{-(k+1)H}A^{-(k+1)})^{-1}]\\
      &\geq \sup_{k} E[(\sum_{i=1}^{\infty}\gamma_i A^{-iH}C^HCA^{-i}+A^{-(k+1)H}A^{-(k+1)})^{-1}] \\
      &= \lim_{k\rightarrow \infty} E[(\sum_{i=1}^{\infty}\gamma_i A^{-iH}C^HCA^{-i}+A^{-(k+1)H}A^{-(k+1)})^{-1}].
    \end{split}
  \end{displaymath}
  Since $A^{-kH}A^{-k} = diag (|\lambda_1|^{-2k},\ldots,|\lambda_n|^{-2k})$, $A^{-kH}A^{-k}$ is monotonically decreasing with respect to $k$ in positive definite sense. Therefore $(\sum_{i=1}^{\infty}\gamma_i A^{-iH}C^HCA^{-i}+A^{-(k+1)H}A^{-(k+1)})^{-1}$ is monotonically increasing. By the Monotone Convergence Theorem, we know that 
  \begin{equation}
    \begin{split}
      \sup_k E[(\widetilde T_k^H \widetilde T_k)^{-1}]& \geq \lim_{k\rightarrow \infty} E[(\sum_{i=1}^{\infty}\gamma_i A^{-iH}C^HCA^{-i}+A^{-(k+1)H}A^{-(k+1)})^{-1}]\\
      & = E[\lim_{k\rightarrow \infty} (\sum_{i=1}^{\infty}\gamma_i A^{-iH}C^HCA^{-i}+A^{-(k+1)H}A^{-(k+1)})^{-1}]\\
      & = E[ (\sum_{i=1}^{\infty}\gamma_i A^{-iH}C^HCA^{-i})^{-1}],
    \end{split}
  \end{equation}
  which proves the left hand side of the inequality.
\end{IEEEproof}

\subsection{Proof of Theorem~\ref{theorem:bound}}

In this part, we will manipulate $(\sum_{i=1}^{\infty}\gamma_i A^{-iH}C^HCA^{-i})^{-1}$ to derive the critical value. The key idea is to use cofactors to find an upper bound of matrix inversion. For non-degenerate system, this upper bound will lead to a upper bound of critical value, which coincides with the lower bound in \cite{Sinopoli:04}. 

Before we prove Theorem~\ref{theorem:bound}, we need the following lemmas.
\begin{lemma}
  \label{lemma:gaussianelimination}
  For a non-degenerate system, it is possible to find a set of row vectors $L_1,L_2,\ldots,L_n$, such that $L_i C = [l_{i,0},\ldots,l_{i,n}]$, where $l_{i,i} = 1$ and $l_{i,j}=0$ if $|\lambda_i|=|\lambda_j|$.
\end{lemma}
\begin{IEEEproof}
  It is simple to show that the lemma holds by using Gaussian Elimination for every equi-block.
\end{IEEEproof}
\begin{lemma}
  \label{lemma:determinant}
  Consider that $|\lambda_{1}| \geq |\lambda_{2}|\cdots \geq |\lambda_n|$, $l_{i,i} = 1$ and $l_{i,j} = 0$ if $i\neq j$ and $|\lambda_i|=|\lambda_j|$. Define $\Delta i_1=i_1$ and $\Delta i_j=i_j-i_{j-1},j=2,3,\ldots,n$. Then the determinant
  \[
  D=\left| {\begin{array}{*{20}c}
    {l_{1,1}\lambda _1^{ - i_1 } } & {l_{1,2}\lambda _2^{ - i_1 } } & \cdots & {l_{1,n}\lambda _n^{ - i_1 } } \\
    {l_{2,1}\lambda _1^{ - i_2 } } & {l_{2,2}\lambda _2^{ - i_2 } } & \cdots & {l_{2,n}\lambda _n^{ - i_2 } } \\
    \vdots & \vdots & \ddots & \vdots \\
    {l_{n,1}\lambda _1^{ - i_n } } & {l_{n,2}\lambda _2^{ - i_n } } & \cdots & {l_{n,n}\lambda _n^{ - i_n } } \\
  \end{array}} \right|
  \]
  is asymptotic to $\prod_{k=1}^n\lambda_k^{-i_k}$, i.e.
  \begin{equation}
    \label{eq:asymptotic}
    \lim_{\Delta i_1,\Delta i_2,\ldots,\Delta i_n \rightarrow \infty} \frac{D}{\prod_{k=1}^n\lambda_k^{-i_k}}=1.
  \end{equation}
\end{lemma}

\begin{lemma}
  \label{lemma:asymptotics}
  If a non-degenerate system satisfies $(H1)-(H3)$ and all its eigenvalues are unstable, then following inequality holds
  \begin{equation}
    \limsup_{\Delta i_1,\ldots,\Delta i_n \rightarrow \infty}\frac{\left(\sum_{j=1}^n A^{-i_jH}C^HCA^{-i_j}\right)^{-1}}{\prod_{j=1}^n |\lambda_j|^{2\Delta i_j}} \leq \overline \beta I_n,
  \end{equation}
  where $\overline \beta > 0$ is a constant, $i_1 < i_2 <\ldots<i_n \in \mathbb N,\;\Delta i_1 = i_1, \Delta i_{j} = i_{j} - i_{j-1}, j= 2,\ldots,n$.
\end{lemma}
\begin{IEEEproof}[Proof of Lemma~\ref{lemma:determinant}]
  The determinant $D$ has $n!$ terms, which have the form $sgn(\sigma) \prod_{k=1}^n l_{k,j_k} \lambda_{j_k}^{-i_k}$. $\sigma = (j_1,j_2,\ldots,j_n)$ is a permutation of set $\{1,2,\ldots,n\}$ and $sgn(\sigma) = \pm 1$ is the sign of permutation. Rewrite \eqref{eq:asymptotic} as
  \[
  \begin{split}
    \frac{D}{\prod_{k=1}^n\lambda_k^{-i_k}}&=\sum_{\sigma} sgn(\sigma) {\frac {\prod_{k=1}^n l_{k,j_k}
    \lambda_{j_k}^{-i_k}}{\prod_{k=1}^n\lambda_k^{-i_k}}} =\sum_{\sigma} sgn(\sigma){\prod_{k=1}^n l_{k,j_k}\frac{\left(\prod_{k=1}^n \lambda_{j_k}\right)^{-\Delta i_1} \cdots\left(\prod_{k=n}^n \lambda_{j_k}\right)^{-\Delta i_n}}{\left(\prod_{k=1}^n \lambda_k\right)^{-\Delta i_1}\cdots\left(\prod_{k=n}^n \lambda_k\right)^{-\Delta i_n}}}\\
    &=\sum_{\sigma} sgn(\sigma){\prod_{k=1}^n l_{k,j_k}\prod_{m=1}^n\left(\frac{\prod_{k=m}^n \lambda_{j_k}}{\prod_{k=m}^n \lambda_k}\right)^{-\Delta i_m}}.
  \end {split}
  \]

  Now we can just analyze each term of the summation. Since $|\lambda_1| \geq \cdots \geq |\lambda_n|$, $|\prod_{k=m}^n \lambda_{j_k}| \geq |\prod_{k=m}^n \lambda_{k}|$. First consider that there exists some $k$s such that $|\lambda_{j_k}| \neq |\lambda_{k}|$ and define $k^*$ to be the largest, which means $|\lambda_{j_{k^*}}| \neq |\lambda_{k^*}|$ and $|\lambda_{j_k}| = |\lambda_{k}|$ for all $k$ greater than $k^*$. Since $|\lambda_{k^*}|$ is the smallest among $|\lambda_1|,\ldots,|\lambda_k|$, we know that $|\lambda_{j_{k^*}}| > |\lambda_{k^*}|$. Thus,
  \begin{displaymath}
    \left| \frac{\prod_{k=k^*}^n \lambda_{j_k}}{\prod_{k=k^*}^n \lambda_k} \right| > 1,
  \end{displaymath}
  and
  \begin{displaymath}
    \lim_{\Delta i_1,\Delta i_2,\ldots,\Delta i_n \rightarrow \infty} \left|\prod_{k=1}^n l_{k,j_k}\prod_{m=1}^n\left(\frac{\prod_{k=m}^n \lambda_{j_k}}{\prod_{k=m}^n \lambda_k}\right)^{-\Delta i_m}\right| \leq |\prod_{k=1}^n l_{k,j_k}| \lim_{\Delta i_{k^*} \rightarrow \infty} \left|\frac{\prod_{k=k^*}^n \lambda_{j_k}}{\prod_{k=k^*}^n \lambda_k}\right|^{-\Delta i_{k^*}} = 0.
  \end{displaymath}

  Then consider that if for all $k$, $|\lambda_{j_k}|=|\lambda_k|$, but $(j_1,\ldots,j_n) \neq (1,2,\ldots,n)$. Thus, there exists $k^*$ such that $j_{k^*} \neq k^*$. Hence $l_{k^*,j_{k^*}} = 0$. Therefore, these terms are always 0.

  The only term left is 
  \begin{displaymath}
    sgn(\sigma)\prod_{k=1}^n l_{k,k}\prod_{m=1}^n\left(\frac{\prod_{k=m}^n \lambda_{k}}{\prod_{k=m}^n \lambda_k}\right)^{-\Delta i_m} = 1.
  \end{displaymath}

  Thus, we can conclude that
  \[
  \lim_{\Delta i_1,\Delta i_2,\ldots,\Delta i_n \rightarrow \infty} \frac{D}{\prod_{k=1}^n\lambda_k^{-i_k}}=1.
  \]
\end{IEEEproof}

\begin{IEEEproof}[Proof of Lemma~\ref{lemma:asymptotics}]
  Because the system is non-degenerate, by Lemma~\ref{lemma:gaussianelimination}, we know that there exist $L_1,L_2,\cdots,L_n$, such that $L_i C = [l_{i,1},\ldots,l_{i,n}]$ is a row vector, $l_{i,i}=1$ and $l_{i,j} = 0$ if $i \neq j$ and $|\lambda_i| =|\lambda_j|$. 

  Define matrices
  \begin{equation}
    U =  \left[ {\begin{array}{*{20}c}
      {l_{1,1}\lambda _1^{ - i_1 } } & {l_{1,2}\lambda _2^{ - i_1 } } & \cdots & {l_{1,n}\lambda _n^{ - i_1 } } \\
      {l_{2,1}\lambda _1^{ - i_2 } } & {l_{2,2}\lambda _2^{ - i_2 } } & \cdots & {l_{2,n}\lambda _n^{ - i_2 } } \\
      \vdots & \vdots & \ddots & \vdots \\
      {l_{n,1}\lambda _1^{ - i_n } } & {l_{n,2}\lambda _2^{ - i_n } } & \cdots & {l_{n,n}\lambda _n^{ - i_n } } \\
    \end{array}} \right], O = U^{-1}.
  \end{equation}

  Define $\overline n_l = \max (\lambda_{max}(L_1^H L_1),\ldots,\lambda_{max}(L_n^H L_n))$. Thus, $L_i^H L_i \leq \overline n_l I_m$, and
  \begin{equation}
    \begin{split}
      &\sum_{j =1}^{n} A^{-i_j H}C^H C A^{-i_j} \geq \sum_{j =1}^{n} \frac{1}{\overline n_l}A^{-i_j H}C^HL_j^HL_j C A^{-i_j}\\
      &=\frac{1}{\overline n_l} \left[ {\begin{array}{*{20}c}
	{A^{-i_1H}C^HL_1^H} & {\cdots} & {A^{-i_nH}C^HL_1^H}  \\
      \end{array}} \right]\left[ {\begin{array}{*{20}c}
	{L_1CA^{-i_1}}  \\
	{\vdots}  \\
	{L_nCA^{-i_n}}  \\
      \end{array}} \right]=\frac{1}{\overline n_l} U^H U,
    \end{split}
  \end{equation}
  and
  \begin{equation}
    \begin{split}
      \label{eq:traceeqsum}
      &\left(\sum_{j =1}^{n} A^{-i_j H}C^H C A^{-i_j}\right)^{-1} \leq \overline n_l \left(U^HU\right)^{-1} = \overline n_l OO^H \leq \overline n_ltrace(OO^H) I_n \\&= \overline n_l\sum_{i,j} O_{i,j}(O^H)_{j,i}I_n = \overline n_l\sum _{i,j} O_{i,j} \times conj(O_{i,j})I_n = \overline n_l\sum_{i,j} |O_{i,j}|^2 I_n,
    \end{split}
  \end{equation}
  where $conj()$ means complex conjugation. 

  Now by Lemma~\ref{lemma:determinant}, we can compute the cofactor matrix of $U$ and hence $O=U^{-1}$. Define the minor $M_{i,j}$ of $U$ as the $(n-1)\times (n-1)$ matrix that results from deleting row $i$ and column $j$. Thus
  \begin{equation}
    O_{i,j}=\frac{(-1)^{i+j}\det(M_{j,i})}{\det (U)}.
    \label{eq:cofactor}
  \end{equation}

  By Lemma~\ref{lemma:determinant}, we know that
  \begin{displaymath}
    \lim_{\Delta i_1,\Delta i_2,\ldots,\Delta i_n \rightarrow \infty} \frac{\det(U)}{\prod_{k=1}^n\lambda_k^{-i_k}}=1.
  \end{displaymath}
  Since $M_{i,j}$ has the same structure with $U$, it is easy to show that 
  \[\det(M_{i,j}) \leq \rho_{i,j}\prod_{k=2}^n |\lambda_k^{-i_{k-1}}|,\] 
  where $\rho_{i,j}$ is a constant. Thus, 
  \begin{equation}
    \begin{split}
      &\limsup_{\Delta i_1,\ldots,\Delta i_n \rightarrow \infty}\frac{\left(\sum_{j =1}^{n} A^{-i_j H}C^H C A^{-i_j}\right)^{-1}}{\prod_{k=1}^n |\lambda_k|^{2\Delta i_k}} \leq \limsup_{\Delta i_1,\ldots,\Delta i_n \rightarrow \infty}\frac{\overline n_l\sum_{i,j}|O_{i,j}|^2}{\prod_{k=1}^n |\lambda_k|^{2\Delta i_k}}I_n\\ 
      &= \limsup_{\Delta i_1,\ldots,\Delta i_n \rightarrow \infty} \overline n_l \left(\sum_{i,j} \left|\frac{\det(M_{i,j})}{\det(U)}\right|^2/\prod_{k=1}^n |\lambda_k|^{2\Delta i_k}\right) I_n\\
      &\leq \overline n_l \left(\sum_{i,j} \rho_{i,j}^2 \left|\frac{\prod_{k=2}^n |\lambda_k^{-i_{k-1}}|}{\prod_{k=1}^n\lambda_k^{-i_k} }\right|^2/\prod_{k=1}^n |\lambda_k|^{2\Delta i_k}\right) I_n = \overline n_l \sum_{i,j}\rho_{i,j}^2 I_n. 
    \end{split}
  \end{equation}

  By \eqref{eq:traceeqsum}, we can conclude that
  \begin{equation}
    \overline \beta  = \overline n_l \sum_{i,j} \rho_{i,j}^2
  \end{equation}
\end{IEEEproof}

Now we are ready to prove Theorem~\ref{theorem:bound}.
\begin{IEEEproof}[Proof of Theorem~\ref{theorem:bound}]
  By Lemma~\ref{lemma:asymptotics}, we know that there exists $\xi_1,\ldots,\xi_n > 0$, such that if $\Delta i_j \geq \xi_j,\; j =1,\ldots,n$, then $ \left[\sum_{j=1}^n A^{-i_jH}C^HCA^{-i_j}\right]^{-1} \leq 2 \overline \beta \prod_{j=1}^n |\lambda_j|^{2\Delta i_j}$. Define stopping time $i_1 = \inf\{i \geq \xi_1|\gamma_i =1\},\;i_j = \inf\{i \geq \xi_j + i_{j-1}|\gamma_i = 1\},\;j=2,\ldots,n$. By the definition of $i_j$ and Lemma~\ref{lemma:asymptotics}, it is simple to show
  \begin{equation}
    \left(\sum_{i=1}^{\infty}\gamma_i A^{-iH}C^HCA^{-i}\right)^{-1} \geq \left(\sum_{j =1}^{n} A^{-i_j H}C^H C A^{-i_j} \right)^{-1}\geq 2\overline \beta \prod_{j=1}^n |\lambda_j|^{2\Delta i_j}.
    \label{eq:finitesummation}
  \end{equation}

  Therefore, $E \left[\left(\sum_{i=1}^{\infty}\gamma_i A^{-iH}C^HCA^{-i}\right)^{-1}\right]$ is bounded if $E \prod_{k=1}^n |\lambda_k^{2\Delta i_{k}}|$ is bounded. From the definition of random variable $i_j$, we know that $\Delta i_j$ are independent of each other. And
  \begin{equation}
    P(\Delta i_j = k) = P (\gamma_{i_{j-1}+\xi_j} = 0,\ldots,\gamma_{i_{j-1}+\xi_j+k-1} = 0,\gamma_{i_{j-1}+\xi_j+k} = 1) = (1-p)^{k-\xi_j-1}p,\quad k \geq \xi_i.
    \label{eq:probability}
  \end{equation}
  Now we can compute the expectation
  \begin{equation}
    \begin{split}
      E \prod_{j=1}^n |\lambda_j^{2\Delta i_{j}}| &= \prod_{j=1}^n E |\lambda_j|^{2\Delta i_j}=\prod_{j=1}^n \sum_{k=\xi_j}^\infty |\lambda_j|^{2 k} P(\Delta i_j = k)\\
      &=\prod_{j=1}^n \sum_{k=\xi_j}^\infty |\lambda_j|^{2 k} (1-p)^{k-\xi_j-1}p,
    \end{split}
  \end{equation}
  which is bounded if and only if 
  \begin{displaymath}
    |\lambda_j|^2 (1-p) < 1 ,\quad j=1,2,\ldots,n.
  \end{displaymath}
  We immediately know that the upper bound for the critical value is $1-|\lambda_1|^{-2}$. Combining with the lower bound given in \cite{Sinopoli:04}, we can complete the proof.
\end{IEEEproof}

Before we finish this subsection, we want to state the following corollaries about the estimation error covariance matrix and boundedness of higher moment of $P_k$:
\begin{corollary}
	\label{corollary:generalineq}
  If a non-degenerate system satisfies $(H1)-(H3)$ and all its eigenvalues are unstable, then the estimation error of state $x_k$ by using only observations $y_{k-i_1} , y_{k - i_2},\ldots,y_{k-i_n}$, where $0 \leq i_1 <  \ldots < i_n \leq k$ and $\Delta i_1 = i_1, \Delta i_{j} = i_{j} - i_{j-1}, j= 2,\ldots,n$, is bounded by
  \begin{equation}
    Cov(x_k|y_{k-i_1} , y_{k - i_2},\ldots,y_{k-i_n}) \leq \overline \beta' \prod_{k=1}^n |\lambda_k|^{2\Delta i_k}I_n,
  \end{equation}
  where $\overline \beta'$ is a constant, provided that $\Delta i_j$ are large enough.
\end{corollary}

\begin{corollary}
	\label{corollary:highmoments}
  If a non-degenerate system satisfies $(H1)-(H3)$ and all its eigenvalues are unstable, then  
  \[
  \sup_k EP_k^{q} \leq \infty \Leftrightarrow p_c>1-|\lambda_1|^{-2q},
  \]
  where $q\in \mathbb N$. 
\end{corollary}

\begin{remark}
	Note that in Corollary~\ref{corollary:generalineq}, we do not assume any distribution of $i_1,\ldots,i_n$. Hence, this corollary allows us to take into account other packet drop models.
\end{remark}
\subsection{Proof of Theorem~\ref{theorem:relaxation} and \ref{theorem:critnondeg}}

In this subsection, we will proof Theorem~\ref{theorem:relaxation} and finally Theorem~\ref{theorem:critnondeg}.

\begin{IEEEproof}[Proof of Theorem~\ref{theorem:relaxation}]
  We consider the case where $R,Q,\Sigma_0$ are identity matrices. To prove the inequality, we will first show that $f_c(A,C) \leq f_c(diag(A_1,A_2),[C_1,C_2])$. Rewrite the system equations as
  \[
  \begin{split}
    \left[ {\begin{array}{*{20}c}
      {x_{k+1,1}}  \\
      {x_{k+1,2}}  \\
    \end{array}} \right] &= \left[ {\begin{array}{*{20}c}
      {A_1} & {}  \\
      {} & {A_2}  \\
    \end{array}} \right]\left[ {\begin{array}{*{20}c}
      {x_{k,1}}  \\
      {x_{k,2}}  \\
    \end{array}} \right] + \left[ {\begin{array}{*{20}c}
      w_{k,1}  \\
      w_{k,2}  \\
    \end{array}} \right] ,\\ 
    x_{k+1,3} &= A_3x_{k,3} + w_{k,3}, \\ 
    y_k &= \left[ {\begin{array}{*{20}c}
      {C_1} & {C_2}  \\
    \end{array}} \right]\left[ {\begin{array}{*{20}c}
      {x_{k,1}}  \\
      {x_{k,2}}  \\
    \end{array}} \right] + v_k + C_3 x_{k,3}. \\ 
  \end{split}
  \]

  Now we want to build a linear filter. Since Kalman filter is the optimal linear filter, the critical value of Kalman filter should be no greater than our linear filter.

  Because $A_3$ is stable, we can just use $\hat x_{k,3} = A_3^k \overline x_{k,3}$ as an unbiased estimation of $x_{k,3}$ and the estimation error covariance is bounded. Now $x_{k,3},x_{k-1,3},\ldots,x_{0,3}$ become measurement noise and we know that
  \begin{equation}
    \begin {split}
    Cov \left(\left[ {\begin{array}{*{20}c}
      {x_{k,3} } \\
      {x_{k - 1,3} } \\
      \vdots \\
      {x_{0,3} } \\
    \end{array}} \right] \right)&= \left[ {\begin{array}{*{20}c}
      I & \cdots & {A_2^{k - 1} } & {A_2^k } \\
      \vdots & \ddots & \vdots & \vdots \\
      0 & \cdots & I & A_2 \\
      0 & \cdots & 0 & I \\
    \end{array}} \right]Cov \left(\left[ {\begin{array}{*{20}c}
      {w_{k - 1,2} } \\
      \vdots \\
      {w_{0,2} } \\
      {x_{0,2} } \\
    \end{array}} \right] \right)\left[ {\begin{array}{*{20}c}
      I & \cdots & {A_2^{k - 1} } & {A_2^k } \\
      \vdots & \ddots & \vdots & \vdots \\
      0 & \cdots & I & A_2 \\
      0 & \cdots & 0 & I \\
    \end{array}} \right]^H.\\
  \end{split}
\end{equation}
Using the same method in Lemma~\ref{lemma:bound} and Theorem~\ref{lemma:pseudoinv} , we can show that this covariance matrix is bounded by $\rho I_{(k+1)n}$, where $\rho$ is a constant independent of $k$. Thus, it is possible to find an i.i.d. sequence of Gaussian measurement noise $v_k'$ such that 
\begin{displaymath}
  Cov \left(\left[ {\begin{array}{*{20}c}
    {v_{k}' } \\
    {v_{k - 1}' } \\
    \vdots \\
    {v_{0}' } \\
  \end{array}} \right] \right) \geq Cov \left(\left[ {\begin{array}{*{20}c}
    {C_3 x_{k,3} } \\
    {C_3 x_{k - 1,3} } \\
    \vdots \\
    {C_3 x_{0,3} } \\
  \end{array}} \right] \right).
\end{displaymath}
We can build another system, 
\begin{equation} 
  \begin{split}
    \left[ {\begin{array}{*{20}c}
      {x_{k+1,1}}  \\
      {x_{k+1,2}}  \\
    \end{array}} \right] &= \left[ {\begin{array}{*{20}c}
      {A_1} & {}  \\
      {} & {A_2}  \\
    \end{array}} \right]\left[ {\begin{array}{*{20}c}
      {x_{k,1}}  \\
      {x_{k,2}}  \\
    \end{array}} \right] + \left[ {\begin{array}{*{20}c}
      w_{k,1}  \\
      w_{k,2}  \\
    \end{array}} \right] ,\\ 
    y_k &= \left[ {\begin{array}{*{20}c}
      {C_1} & {C_2}  \\
    \end{array}} \right]\left[ {\begin{array}{*{20}c}
      {x_{k,1}}  \\
      {x_{k,2}}  \\
    \end{array}} \right] + v_k + v_k' .\\ 
  \end{split}
  \label{eq:systemdiscription3}
\end{equation}
By the property of linear filter, the estimation error of Kalman filter for system \eqref{eq:systemdiscription3} will be greater than the one for the original system, which implies that
\begin{equation}
  f_c(A,C) \leq f_c(diag(A_1,A_2),[C_1,C_2]).  
  \label{eq:inequality1}
\end{equation}

Now define function
\[
g(X,A,\gamma)=A X A^H+Q-\gamma A X C (C X C^H+R)^{-1}C X A^H,
\]
and $\alpha,\beta$ to be scalars. Therefore 
\[
g(X,\alpha A,\gamma)-g(X,\beta A,\gamma)=(\alpha^2-\beta^2)[AXA^H-\gamma AXC^H(CXC^H+R)^{-1}CXA^H].
\]

Thus, $g(X,\alpha A,\gamma)$ is a non-decreasing function of $\alpha$ when $\alpha > 0$. Since for system $(A, C , R, Q, \Sigma_0)$, the error covariance matrix $P_k$ follows recursive equation $P_{k+1} = g(P_{k}, A, \gamma_k)$. By the monotonicity of $g(X,\alpha A,\gamma_k)$, we know that $P_{k}$ is also a non-decreasing function of $\alpha$. Thus, the critical value for the system is also non-decreasing, which implies that
\begin{equation}
  f_c(A,C) \leq  \lim_{\alpha \rightarrow 1^+} f_c(\alpha A, C).
  \label{eq:inequality2}
\end{equation}
The limit on the right hand side must exist because of the monotonicity of function $f_c$. Combining \eqref{eq:inequality1} and \eqref{eq:inequality2}, we can finish the proof. 
\end{IEEEproof}

\begin{IEEEproof}[Proof of Theorem~\ref{theorem:critnondeg}]
  If the system does not have unstable and critically stable eigenvalues, then the proof is trivial. Otherwise by Theorem~\ref{theorem:criticalfunction}, we know that 
  \begin{displaymath}
    f_c(A,C) \geq f_c(A_1,C_1) = \max\{1-|\lambda_1|^{-2},0\}.
  \end{displaymath}
  By Theorem~\ref{theorem:relaxation}, 
  \begin{displaymath}
    f_c(A,C) \leq \lim_{\alpha\rightarrow 1^+} f_c(diag(\alpha A_1,\alpha A_2),[C_1\;C_2]) = \lim_{\alpha \rightarrow 1^+}\max\{1-|\alpha\lambda_1|^{-2},0\} = \max\{1-|\lambda_1|^{-2},0\}.
  \end{displaymath}
  Hence, the critical value is $1-|\lambda_1|^{-2}$.
\end{IEEEproof}

\section{A Complete Characterization of Critical Value for Second Order Systems}
\label{sec:2degsys}
This section is devoted to a complete characterization of linear system with a diagonal $2$ by $2$ $A$ matrix. This kind of systems can be seen as the building blocks of larger systems and thanks to Theorem~\ref{theorem:criticalfunction}, we know that the critical value of such system will give a lower bound of critical value for larger system. To our surprise, the critical values of second order systems are in fact quite complex.

For non-degenerate second order systems, we can directly apply Theorem~\ref{theorem:critnondeg} to derive the critical value. Thus, we will focus on degenerate systems. Using the same strategy as the previous section, we will first deal with unstable degenerate system, then apply Theorem~\ref{theorem:relaxation} to generalize the result to critically stable and stable systems. By the definition of degeneracy, we know that a detectable second order system is degenerate if and only if the following assumptions holds:
\begin{enumerate}
  \item $\lambda_2 = \lambda_1 exp(j\varphi) $ , where $j^2 = -1$ and $\varphi \in (0,2\pi)$. ($\varphi \neq 0$, otherwise the system is not detectable)
  \item $rank(C) = 1$.
\end{enumerate}

To simplify the notation, let us define $\lambda \triangleq |\lambda_1| = |\lambda_2|$, $z \triangleq exp(j\varphi)$. The proof of critical value is divided into 2 parts. First we want to deal with the case when $\varphi/2\pi$ is rational: 
\begin{theorem}
  \label{theorem:criticalvaluefordegsystem}
  If a unstable second order degenerate system satisfies hypothesis $(H1) - (H3)$ , then the critical value of the system is 
  \begin{equation}
    p_c = 1-|\lambda_1|^{-2q/q-1},
    \label{eq:critvaluefor2degenerate}
  \end{equation}
  where $\varphi/2\pi = r/q$ , $q > r$ and $r,q \in \mathbb N$ are irreducible.
\end{theorem}
Then we consider the case when $\varphi/2\pi$ is irrational:
\begin{theorem}
  \label{theorem:criticalvaluefordegsystem2}
  If a unstable second order degenerate system satisfies hypothesis $(H1)-(H3)$, then the critical value of the system is 
  \begin{equation}
    p_c = 1-|\lambda_1|^{-2},
    \label{eq:critvaluefor2degenerate2}
  \end{equation}
  if $\varphi/2\pi$ is irrational.
\end{theorem}
\begin{IEEEproof}[Proof of Theorem~\ref{theorem:criticalvaluefordegsystem}]
  By the properties of degeneracy, we know that $rank(C^H C) = 1$. Thus, $C^HC = \left[{\begin{array}{*{20}c}
    a^2 & ab\\
    ab & b^2\\
  \end{array}}\right]$, where $a,b$ are real constants. It can be also proved that $a,b \neq 0$ due to the detectability of $(C,A)$. Since $A = diag(\lambda_1,\lambda_2) = diag(\lambda_1,\lambda_1z) $, we know that
  \begin{equation}
    \begin{split}
      \sum_{i=1}^\infty \gamma_i A^{-iH}C^HCA^{-i} &= \sum_{i=1}^\infty \gamma_i \left[{\begin{array}{*{20}c}
	a^2\lambda^{-2i} & ab \lambda^{-2i} z^{-i}\\
	ab\lambda^{-2i}z^i & b^2 \lambda^{-2i}\\
      \end{array}}\right]\\ &= \left[{\begin{array}{*{20}c}
	a & \\
	& b\\
      \end{array}}\right]\left( \sum_{i=1}^\infty \gamma_i \lambda^{-2i} \left[{\begin{array}{*{20}c}
	1 &  z^{-i}\\
	z^i & 1\\
      \end{array}}\right]\right)  \left[{\begin{array}{*{20}c}
	a & \\
	& b\\
      \end{array}}\right].
    \end{split}
  \end{equation}

  Since $a,b \neq 0$, we know that $E(\sum_{i=1}^\infty \gamma_i A^{-iH}C^HCA^{-i})^{-1}$ is bounded if and only if
  \[E\left(\sum_{i=1}^\infty \gamma_i \lambda^{-2i} \left[{\begin{array}{*{20}c}
    1 &  z^{-i}\\
    z^i & 1\\
  \end{array}}\right]\right)^{-1} < \infty.\]

  Define
  \begin{displaymath}
   \Xi \triangleq \sum_{i=1}^\infty \gamma_i \lambda^{-2i} \left[{\begin{array}{*{20}c}
      1 &  z^{-i}\\
      z^i & 1\\
    \end{array}}\right] .
  \end{displaymath}

  It is easy to show that
  \begin{equation}
    \begin{split}
      \label{eq:traceanddet42degsys}
      trace(\Xi) &= 2\sum_{i=1}^\infty \gamma_i\lambda^{-2i},\\
      \det(\Xi) &= \left(\sum_{i=1}^\infty \gamma_i \lambda^{-2i}\right)^2 - \left(\sum_{i=1}^{\infty} \gamma_i \lambda^{-2i} z^i\right) \times \left(\sum_{i=1}^{\infty} \gamma_i \lambda^{-2i} z^{-i}\right)\\
      &=\sum_{i=1}^{\infty} \gamma_i \lambda^{-4i} + 2\sum_{i=1}^{\infty} \sum_{j=i+1}^{\infty} \gamma_i\gamma_j \lambda^{-2i}\lambda^{-2j} - \sum_{i=1}^{\infty} \gamma_i \lambda^{-4i} -\sum_{i=1}^{\infty} \sum_{j=i+1}^{\infty} \gamma_i\gamma_j \lambda^{-2i}\lambda^{-2j}(z^{i-j}+z^{j-i})\\
      & = \sum_{i=1}^{\infty} \sum_{j=i+1}^{\infty} \gamma_i\gamma_j \lambda^{-2i}\lambda^{-2j}(2-z^{i-j}-z^{j-i}).\\
    \end{split}
  \end{equation}

  Define set $\mathbb S_{q,\infty} = \{l \in \mathbb N|l \neq kq, k \in \mathbb N\}$ and $\mathbb S_{q,i} = \{l \in \mathbb S_{q,\infty}|l < i\}$. Since $z = exp(2r\pi/q)$ and $q,r$ are irreducible, $z^{j-i} = 1$ if and only if $ |j-i| \notin \mathbb S_{q,\infty}$. It is easy to show that 
  \[
  \inf\{2-z^{i-j}-z^{j-i}||j - i| \in \mathbb S_{q,\infty}\} = \inf \{2-z^i-z^{-i}|i = 1,\ldots,q-1\} = 2 - 2\cos(\frac{2\pi}{q}) > 0,
  \]
  and
  \[
  \sup\{2-z^{i-j}-z^{j-i}||j - i| \in \mathbb S_{q,\infty}\} = \sup \{2-z^i-z^{-i}|i = 1,\ldots,q-1\} < 4.
  \]
  Thus, 
  \begin{equation}
    \begin{split}
      \label{eq:boundfordet1}
      [2-2\cos(2\pi/q)] \sum_{i=1}^\infty \sum_{j-i \in \mathbb S_{q,\infty}} \gamma_i\gamma_j\lambda^{-2i}\lambda^{-2j} &\leq\det(\Xi) = \sum_{i=1}^\infty \sum_{j-i \in \mathbb S_{q,\infty}} \gamma_i\gamma_j\lambda^{-2i}\lambda^{-2j} \\
      &\leq 4 \sum_{i=1}^\infty \sum_{j-i \in \mathbb S_{q,\infty} } \gamma_i\gamma_j\lambda^{-2i}\lambda^{-2j}.
    \end{split}
  \end{equation}

  Define stopping time $\tau_1 = \inf\{i \in \mathbb N|\gamma_i = 1\}$ and $\tau_2 = \inf \{j\in \mathbb N| j - \tau_1 \in \mathbb S_{q,\infty},\gamma_j = 1\}$. Thus,
  \begin{equation}
    \label{eq:boundfortrace}
    \lambda^{-2\tau_1} \leq trace(\Xi) = \sum_{i=1}^{\infty}\gamma_i\lambda^{-2i} = \sum_{i = \tau_1}^{\infty} \gamma_i\lambda^{-2i} \leq \sum_{i=\tau_1}^{\infty} \lambda^{-2i} = \frac{1}{1-\lambda^{-2}}\lambda^{-2\tau_1}.
  \end{equation}

  Now consider there exist two index $a,b$ such that $b > a, b-a \in \mathbb S_{q,\infty}$ and $\gamma_a = \gamma_b = 1$. By the definition of $\tau_1$, we know that $\tau_1 \leq a$. Suppose that $b < \tau_2$, therefore $\tau_1 \leq a < b < \tau_2$. By the definition of $\tau_2$, $a-\tau_1 = k_a q ,b - \tau_1 = k_b q$. As a result, $b-a = (k_b-k_a)q$, which contradicts with the fact $b-a\in \mathbb S_{q,\infty}$. Therefore we can conclude that $\tau_2 \leq b$. 

  Thus, for all $\gamma_a\gamma_b = 1, b-a \in \mathbb S_{q,\infty}$,$\tau_1 \leq a,\tau_2 \leq b$, which gives
  \begin{equation}
    \begin{split}
      \label{eq:boundfordet2}
      \lambda^{-2\tau_1}\lambda^{-2\tau_2} &\leq \sum_{i=1}^\infty \sum_{j-i \in \mathbb S_{q,\infty}} \gamma_i\gamma_j\lambda^{-2i}\lambda^{-2j} = \sum_{i = \tau_1}^\infty \sum_{j \geq \tau_2,j - i \in \mathbb S{q,\infty}} \gamma_i\gamma_j\lambda^{-2i}\lambda^{-2j}\\
      & \leq \sum_{i=\tau_1}^{\infty}\sum_{j=\tau_2}^{\infty} \lambda^{-2i}\lambda^{-2j}= \frac{1}{(1-\lambda^{-2})^2}\lambda^{-2\tau_1}\lambda^{-2\tau_2}.
    \end{split}
  \end{equation}

  Define $\sigma_1,\sigma_2$ to be the eigenvalues of $\Xi$. Thus,
  \begin{equation}
    \label{eq:traceeqtraceoverdet}
    trace(\Xi^{-1}) = \sigma_1^{-1} + \sigma_2^{-1} = \frac{\sigma_1 + \sigma_2}{\sigma_1\sigma_2} = \frac{trace(\Xi)}{\det(\Xi)}.
  \end{equation}

  By inequality \eqref{eq:boundfordet1}, \eqref{eq:boundfortrace} and \eqref{eq:boundfordet2}, it is easy to justify that $E trace(\Xi^{-1}) < \infty$ if and only if 
  \begin{displaymath}
    E\frac{\lambda^{-2\tau_1}}{\lambda^{-2\tau_1}\lambda^{-2\tau_2}} = E\lambda^{2\tau_2}  =  E(\lambda^{2\tau_1}\lambda^{2(\tau_2-\tau_1)})< \infty.
  \end{displaymath}

  Now we need to compute the distribution of $\tau_1, \tau_2 - \tau_1$. By definition, the event $\{\tau_1 = i\}$ is equivalent to $\{\gamma_1 = \ldots = \gamma_{i-1} = 0, \gamma_i = 1\}$. The event $\{\tau_2 - \tau_1 = i\}$, where $i \in \mathbb S_{q,\infty}$, is equivalent to $\{\gamma_{\tau_1 + j} = 0, \gamma_{\tau_2} = 1 \}$, for all $j \in \mathbb S_{q,i}$. Since $\tau_2 - \tau_1$ only depends on $\gamma_{\tau_1+i}, i \in \mathbb S_{q,\infty}$, $\tau_2 - \tau_1$ is independent of $\tau_1$. The distributions of $\tau_1,\tau_2 - \tau_1$ are characterized by the following equations:
  \begin{equation}
    P(\tau_1 = i) = P(\gamma_1=\ldots=\gamma_{i-1} = 0, \gamma_i = 1) = (1-p)^{i-1}p,
  \end{equation}
  and
  \begin{equation}
    P(\tau_2 - \tau_1 = i) = P(\gamma_{\tau_1 + j} = 0,\gamma_{\tau_2} = 1) = (1-p)^{|\mathbb S_{q,i}|}p,
  \end{equation}
  where $j \in \mathbb S_{q,i}, i \in \mathbb S_{q,\infty}$ and $|\mathbb S_{q,i}|$ means the number of elements in $\mathbb S_{q,i}$. Thus,
  \begin{equation}
    E\left(\lambda^{2\tau_1}\lambda^{2(\tau_2-\tau_1)}\right)=E\lambda^{2\tau_1} \times E\lambda^{2(\tau_2-\tau_1)} = \sum_{i = 1}^{\infty}p (1-p)^i\lambda^{2i} \times \sum_{i \in \mathbb S_{q,\infty}} p(1-p)^{|\mathbb S_{q,i}|}\lambda^{2i}.
    \label{eq:seriesproductfordegsystem}
  \end{equation}

  The first series is a simple geometric series which is bounded if and only if $p > 1-1/\lambda^2$. Using the root test of convergence, $\sum_{i \in \mathbb S_{q,\infty}} p(1-p)^{|\mathbb S_{q,i}|}\lambda^{2i}$ is bounded if and only if 
  \begin{equation}
    \limsup_{|\mathbb S_{q,i}| \rightarrow \infty} \sqrt[|\mathbb S_{q,i}|]{p(1-p)^{|\mathbb S_{q,i}|}\lambda^{2i}} = (1-p)\limsup_{|\mathbb S_{q,i}| \rightarrow \infty}\lambda^{2\frac{i}{|\mathbb S_{q,i}|}} < 1.
  \end{equation}

  Since $|\mathbb S_{q,i}| = \lceil (i-1)(q-1)/q \rceil$, where $\lceil x \rceil$ means the minimal integer that is no less that $x$, $\limsup_{|\mathbb S_{q,i} \rightarrow \infty|} \left(i/|\mathbb S_{q,i}|\right) = q/(q-1)$. As a result, the second series convergences if and only if 
  \begin{equation}
    (1-p)\lambda^{2q/(q-1)} < 1,
  \end{equation}
  which is equivalent to $p > 1-\lambda^{-2q/(q-1)}$. Now we can conclude that the critical arrival probability is 
  \begin{equation}
    p_c = 1-\lambda^{-\frac{2q}{q-1}}.
  \end{equation}
\end{IEEEproof}

\begin{IEEEproof}[Proof of Theorem~\ref{theorem:criticalvaluefordegsystem2}]
  The proof is quite similar to the proof of Theorem~\ref{theorem:criticalvaluefordegsystem}. The proof before \eqref{eq:traceanddet42degsys} still holds. However, in \eqref{eq:boundfordet1} we need to change the set $\mathbb S_{q,\infty}$. Define set $\mathbb T_{\varepsilon,\infty} = \{l \in \mathbb N|2 - z^{l} - z^{-l} > \varepsilon\}$, where $\varepsilon > 0$. And $\mathbb T_{\varepsilon,i} = \{l \in \mathbb T_{\varepsilon,\infty}| l < i\}$. Therefore, \eqref{eq:boundfordet1} becomes
  \begin{equation}
    \label{eq:boundfordet3}
    \varepsilon \sum_{i=1}^\infty \sum_{j-i \in \mathbb T_{\varepsilon,\infty}} \gamma_i\gamma_j\lambda^{-2i}\lambda^{-2j} \leq det(\Xi) = \sum_{i=1}^\infty \sum_{j > i} \gamma_i\gamma_j\lambda^{-2i}\lambda^{-2j}(2 - z^{i-j} -z ^{j-i}) .
  \end{equation}
  \eqref{eq:boundfortrace}, \eqref{eq:traceeqtraceoverdet} still hold if we change every set $\mathbb S_{q,i}$ to $\mathbb T_{\varepsilon,i}$. However only the left side inequality in \eqref{eq:boundfordet2} holds, because there is no guarantee that for all $a,b$ satisfies $b - a \in \mathbb T_{\varepsilon,\infty}, \gamma_a = \gamma_b = 1$, $\tau_1 \leq a,\tau_2 \leq b$ always holds. Also in Inequality \eqref{eq:boundfordet3}, we only prove the left side inequality of \eqref{eq:boundfordet1}. As a result, $E\lambda^{2\tau_2} < \infty$ will only be the sufficient condition for the boundedness of estimation error covariance. Following the rest of the proof, it can be derived that
  \begin{equation}
    p \geq 1 - \limsup_{|\mathbb T_{\varepsilon,i}| \rightarrow \infty }\lambda^{-2\frac{i}{|\mathbb T_{\varepsilon,i}|}}
  \end{equation}
  is sufficient for bounded estimation error. Since $\varepsilon$ can be any positive real number, we can conclude that
  \begin{equation}
    \label{eq:upperboundfor2degsys}
    p_c \leq 1- \lim_{\varepsilon \rightarrow 0^+} \limsup_{|\mathbb T_{\varepsilon,i}|\rightarrow \infty} \lambda^{-2\frac{i}{|\mathbb T_{\varepsilon,i}|}}.
  \end{equation}

  Now we only need to estimate $i/|\mathbb T_{\varepsilon,i}|$. $2 - z^i - z^{-i} = 2 - 2\cos(i\varphi)$. Thus 
  \[
  \begin{split}
    2 - 2 \cos (i\varphi) \geq \varepsilon &\Leftrightarrow i \varphi \notin [2k\pi - \arccos(1-\varepsilon/2),2k\pi+\arccos(1-\varepsilon/2)], k \in \mathbb Z\\
    &\Leftrightarrow i(\varphi/2\pi) \notin [k - \Delta_\varepsilon, k+\Delta_\varepsilon], k \in \mathbb Z,
  \end{split}
  \]
  where $\Delta_\varepsilon = \arccos(1-\varepsilon/2)/2\pi$.

  Define $N_\varepsilon = \inf \{i\in \mathbb N| i(\varphi/2\pi) \in [k - 2\Delta_\varepsilon, k + 2\Delta_\varepsilon], k \in \mathbb Z\}$. Suppose that $a(\varphi/2\pi),b(\varphi/2\pi),b > a,$ both belong to $[k-\Delta_\varepsilon,k+\Delta_\varepsilon],k\in \mathbb Z$. Thus, $(b-a)(\varphi/2\pi) \in [k-2\Delta_\varepsilon,k+2\Delta_\varepsilon], k \in \mathbb Z$. By the definition of $N_\varepsilon$, we can conclude that $b - a \geq N_\varepsilon$, which implies that if $a(\varphi/2\pi) \in [k-\Delta_\varepsilon,k+\Delta_\varepsilon]$, then $(a+1)(\varphi/2\pi),\ldots, (a+N_\varepsilon -1)(\varphi/2\pi) \notin [k-\Delta_\varepsilon,k+\Delta_\varepsilon]$. Therefore, if $a \notin \mathbb T_{\varepsilon,\infty}$, then $a+1,\ldots,a+N_\varepsilon-1 \in \mathbb T_{\varepsilon,\infty}$. As a result,
  \begin{equation}
    \frac{N_\varepsilon}{N_\varepsilon - 1} \geq \limsup_{|\mathbb T_{\varepsilon,i}| \rightarrow \infty} \frac{i}{|\mathbb T_{\varepsilon,i}|} \geq 1.
  \end{equation}

  Since $\varphi/2\pi$ is irrational, $\lim_{\varepsilon\rightarrow 0^+} N_\varepsilon = \infty$. Therefore,
  \begin{equation}
    \label{eq:limforirrantionalnumberproximation}
    \lim_{\varepsilon\rightarrow 0^+}\limsup_{|\mathbb T_{\varepsilon,i} \rightarrow \infty|} \frac{i}{|\mathbb T_{\varepsilon,i}|} = 1.
  \end{equation}

  By \eqref{eq:upperboundfor2degsys} and \eqref{eq:limforirrantionalnumberproximation}, we can conclude that the critical arrival probability $p_c$ satisfies
  \begin{displaymath}
    p_c \leq 1-\lambda^{-2},
  \end{displaymath}
  which is exactly the lower bound in \cite{Sinopoli:04}. Therefore, we can conclude the proof.
\end{IEEEproof}

Now we can proof the main theorem:
\begin{IEEEproof}[Proof of Theorem~\ref{theorem:critdeg}]
  By Theorem~\ref{theorem:critnondeg}, \ref{theorem:criticalvaluefordegsystem} and \ref{theorem:criticalvaluefordegsystem2}, we know that the only case we need to prove is critically stable degenerate systems, which is trivial by directly applying Theorem~\ref{theorem:relaxation}.
\end{IEEEproof}

\section{Conclusions and Future Work}
\label{sec:conclusion}

In this paper we address the problem of state estimation for a discrete-time linear Gaussian system where observations are communicated to the estimator via a memoryless erasure channel.  Following the work of Sinopoli et Al.~\cite{Sinopoli:04}, we were able to compute the value of the critical probability for a very general class of linear systems. The boundedness analysis in this paper can be easily generalized to general Markovian packet loss models and to the boundedness of higher moments of the error covariance. Future work will attempt at determining the complete statistics of the error covariance matrix of the Kalman Filter under Bernoulli losses.

\section{ACKNOWLEDGMENTS}
\label{sec:ack}

The authors gratefully acknowledge Professors Francesco Bullo and P. R. Kumar, Craig Robinson for the numerous interesting discussion on the topic.

\bibliographystyle{IEEEtran}
\bibliography{NetContSys}

\begin{thebibliography}{10}
\providecommand{\url}[1]{#1}
\csname url@samestyle\endcsname
\providecommand{\newblock}{\relax}
\providecommand{\bibinfo}[2]{#2}
\providecommand{\BIBentrySTDinterwordspacing}{\spaceskip=0pt\relax}
\providecommand{\BIBentryALTinterwordstretchfactor}{4}
\providecommand{\BIBentryALTinterwordspacing}{\spaceskip=\fontdimen2\font plus
\BIBentryALTinterwordstretchfactor\fontdimen3\font minus
  \fontdimen4\font\relax}
\providecommand{\BIBforeignlanguage}[2]{{%
\expandafter\ifx\csname l@#1\endcsname\relax
\typeout{** WARNING: IEEEtran.bst: No hyphenation pattern has been}%
\typeout{** loaded for the language `#1'. Using the pattern for}%
\typeout{** the default language instead.}%
\else
\language=\csname l@#1\endcsname
\fi
#2}}
\providecommand{\BIBdecl}{\relax}
\BIBdecl

\bibitem{Sinopoli:04}
B.~Sinopoli, L.~Schenato, M.~Franceschetti, K.~Poolla, M.~Jordan, and
  S.~Sastry, ``Kalman filtering with intermittent observations,'' \emph{IEEE
  Transactions on Automatic Control}, vol.~49, no.~9, pp. 1453--1464, September
  2004.

\bibitem{wireless_sensor_network}
N.~P. Mahalik, \emph{Sensor Networks and Configuration}.\hskip 1em plus 0.5em
  minus 0.4em\relax Springer, 2007.

\bibitem{kp-fb:07j}
K.~Plarre and F.~Bullo, ``On kalman filtering for detectable systems with
  intermittent observations,'' \emph{IEEE Transactions on Automatic Control},
  vol.~54, no.~2, pp. 386--390, 2009.

\bibitem{Yilin08}
Y.~Mo and B.~Sinopoli, ``A characterization of the critical value for kalman
  filtering with intermittent observations,'' in \emph{Proceedings of IEEE
  Conference on Decision and Control.}, 2008.

\bibitem{Xu:05}
Y.~Xu and J.~Hespanha, ``Estimation under controlled and uncontrolled
  communications in networked control systems,'' in \emph{Proceedings of the
  CDC-ECC}, Sevilla, Spain, December 2005, pp. 842--847.

\bibitem{Craig:07}
C.~Robinson and P.~R. Kumar, ``Sending the most recent observation is not
  optimal in networked control: Linear temporal coding and towards the design
  of a control specific transport protocol,'' in \emph{Proc. of IEEE CDC}, Dec
  2007, pp. 334--339.

\bibitem{Minyi:07}
M.~Huang and S.~Dey, ``Stability of kalman filtering with markovian packet
  losses,'' \emph{Automatica}, vol.~43, no.~4, pp. 598--607, 2007.

\bibitem{Liu:04}
X.~Liu and A.~Goldsmith, ``Kalman filtering with partial observation losses,''
  in \emph{Proceedings of IEEE Conference on Decision and Control}, vol.~4,
  Bahamas, December 2004, pp. 4180--4186.

\bibitem{xielihua09}
X.~Lihua and L.~Shi, ``State estimation over unreliable network,'' in
  \emph{Asian Control Conference}, August 2009.

\bibitem{censi09}
A.~Censi, ``On the performance of {K}alman filtering with intermittent
  observations: a geometric approach with fractals,'' in \emph{Proceedings of
  the American Control Conference (ACC)}, 2009.

\bibitem{Ali:09}
A.~Vakili and B.~Hassibi, ``On the steady-state performance of kalman filtering
  with intermittent observations for stable systems,'' in \emph{8th IEEE
  Conference on Decision and Control (CDC)}, Dec 2009.

\bibitem{Gupta:05}
V.~Gupta, T.~Chung, B.~Hassibi, and R.~M. Murray, ``On a stochastic sensor
  selection algorithm with applications in sensor scheduling and sensor
  coverage,'' \emph{Automatica}, vol.~42, no.~2, pp. 251--260, 2006.

\bibitem{Ling09}
L.~Shi and L.~Qiu, ``State estimation over a network: Packet-dropping analysis
  and design,'' in \emph{The 7th IEEE International Conference on Control and
  Automation}, Dec 2009.

\bibitem{gerschgorincircle}
R.~S. Varga, \emph{Ger\v sgorin and His Circles}.\hskip 1em plus 0.5em minus
  0.4em\relax New York: Springer, 2004.

\end{thebibliography}
\end{document}